\documentclass[12pt]{article}
\usepackage{amssymb}
\usepackage{amsmath}
\usepackage{amsfonts}
\usepackage{xcolor}

\newtheorem{theorem}{Theorem}[section]

\newtheorem{corollary}{Corollary}[section]

\newtheorem{definition}{Definition}[section]
\newtheorem{example}{Example}[section]

\newtheorem{lemma}{Lemma}[section]

\newtheorem{remark}{Remark}[section]

\newenvironment{proof}[1][Proof]{\noindent \textbf{#1.} }{\  \rule{0.5em}{0.5em}}
\begin{document}

\title{The role of certain Brauer and Rado results in the nonnegative
inverse spectral problems\thanks{
Supported by CONICYT-FONDECYT 1170313, Chile, and CONICYT-PAI 79160002,
Chile.}}
\author{Ana I. Julio\thanks{%
Corresponding author: ajulio@ucn.cl (A. I. Julio), rsoto@ucn.cl (R.L.Soto) }%
, Ricardo L. Soto \\
Departamento de Matem\'{a}ticas\\
Universidad Cat\'{o}lica del Norte \\
Casilla 1280, Antofagasta, Chile.}
\date{}
\maketitle

\begin{abstract}
\ \newline
We say that a list $\Lambda =\{ \lambda _{1},\ldots ,\lambda _{n}\}$ of
complex numbers is realizable, if it is the spectrum of a nonnegative matrix
$A$ (the realizing matrix). We say that $\Lambda $ is universally realizable
if it is realizable for each possible Jordan canonical form allowed by $%
\Lambda .$ This work does not contain new results. As its title says, our
goal is to show and emphasize the relevance of certain results of Brauer and
Rado in the study of nonnegative inverse spectral problems. We show that
virtually all known results, which give sufficient conditions for the list $%
\Lambda $ to be realizable or universally realizable, can be obtained from
the results of Brauer or Rado. Moreover, in this case, we may always compute
a realizing matrix.
\end{abstract}


\bigskip

\section{Introduction}

\noindent The \textit{nonnegative inverse eigenvalue problem} (hereafter
\textit{NIEP}) is the problem of characterizing all possible spectra of
entrywise nonnegative matrices. If there exists a nonnegative matrix $A$
with spectrum $\Lambda =\{\lambda _{1},\lambda _{2},\ldots ,\lambda _{n}\},$
we say that $\Lambda $ is \textit{realizable} and that $A$ is the realizing
matrix. The \textit{NIEP }remains unsolved for $n\geq 5.$ In the general
case, when the possible spectrum $\Lambda $ is a list of complex numbers,
the problem has been solved for $n=3$ by Loewy and London \cite{Loewy}, and
for $n=4$ by Meehan \cite{Meehan}, and independently by Torre-Mayo et al.
\cite{Torre}. The case $n=5$ have been solved for realizing matrices of
trace zero by Laffey and Meehan \cite{Laffey2}. When $\Lambda $ is a list of
real numbers, the \textit{NIEP} is called the \textit{real nonnegative
inverse eigenvalue problem }(\textit{RNIEP}), and a number of sufficient
conditions for the existence of a solution are known (see \cite{Soto, Soto4,
Soto7} and the references therein). \ If the realizing matrix is required to
be symmetric we have the \emph{symmetric nonnegative inverse eigenvalue
problem }(\emph{SNIEP}), which has been solved for $n=5$ with realizing
matrices of trace zero by Spector \cite{Spector}. For $n\leq 4,$ the \emph{%
RNIEP} and the \textit{SNIEP }are equivalent, while for $n\geq 5,$ they are
different \cite{Johnson2}. A number of sufficient conditions for the
existence of a symmetric nonnegative matrix with prescribed spectrum have
also been obtained (see \cite{Soto5, Soto6, Soto7} and the references
therein).\\ \ \\
We say that a list of complex numbers $\Lambda =\{\lambda
_{1},\lambda _{2},\ldots ,\lambda _{n}\}$ is universally realizable (\textit{%
UR}), if $\Lambda $ is realizable for each possible Jordan canonical form (%
\textit{JCF}) allowed by the list $\Lambda .$ The problem of the universal
realizability of spectra will be called the universal realizability problem (%
\textit{URP}).\\ \ \\
A set $\mathcal{K}$ of conditions is said to be a realizability criterion%
\textit{\ }if any list $\Lambda $ of complex numbers satisfying the
conditions $\mathcal{K}$ is realizable. A real matrix \thinspace $%
A=(a_{ij})_{i=1}^{n}$\thinspace\ is said to have \textit{constant row sums}
if all its rows sum up to a same constant, say $\alpha $, that is, $%
\sum_{j=1}^{n}a_{ij}=\alpha ,\ \ i=1,\ldots ,n.$ The set of all real
matrices with constant row sums equal to $\alpha $\ will be denoted by $%
\mathcal{CS}_{\alpha }$. It is clear that any matrix in $\mathcal{CS}%
_{\alpha }$ has an eigenvector $\mathbf{e}=(1,\ldots ,1)^{T}$ corresponding
to the eigenvalue $\alpha $. We denote by $\mathbf{e}_{i}=(0,\ldots
,1,\ldots ,0)^{T},$\ with the $1$ in the $i^{th}$ position, the $i^{th}$
column of the identity matrix of the appropriate size. The importance of the
real matrices with constant row sums is due to the fact that the problem of
finding a nonnegative matrix $A$ with spectrum $\Lambda =\{\lambda
_{1},\lambda _{2},\ldots ,\lambda _{n}\},$ $\lambda _{1}$ being the Perron
eigenvalue, is equivalent to the problem of finding a nonnegative matrix in $%
\mathcal{CS}_{\lambda _{1}}.$\\ \ \\
The purpose of this work is to examine the \textit{NIEP} and the \textit{URP}, from Brauer and Rado results point of view \cite{Brauer, Perfect2}. In
particular, we show that virtually all known results, which give
realizability criteria for the \textit{NIEP} and the \textit{URP} to have a
solution, may be obtained by applying certain results due to Brauer or Rado,
which we identify in Section $2.$ Moreover, the proofs from Brauer and Rado
results are constructive, in the sense that they always allow us to compute
a realizing matrix.\\ \ \\
Brauer's result, Theorem \ref{TeoBrauerAA} below, shows how to modify one
single eigenvalue of a matrix via a rank-one perturbation, without changing
any of the remaining eigenvalues. This, together with the properties of real
matrices with constant row sums, are the basic ingredients of the technique
that have been used in most cases and it suggests that Brauer result,
Theorem \ref{TeoBrauerAA} in Section $2,$ can be a very useful tool to deal
with the \textit{NIEP} \emph{\ }and the \textit{URP}. This approach goes
back to Perfect who first used it in \cite{Perfect1} to obtain sufficient
conditions, but it was somehow abandoned for many years until in \cite{Soto}%
, the author rediscovered it to obtain sufficient conditions for the
realizability of partitioned real spectra, with the partition allowing some
of its pieces to be nonrealizable. \\ \ \\
Rado's result, Theorem \ref{TeoBrauerC} in Section $2$, is an extension of
Theorem \ref{TeoBrauerAA} and shows how to modify $r$ eigenvalues of a
matrix of order $n,$ $r<n,$ via a rank-$r$ perturbation, without changing
any of the remaining $(n-r)$ eigenvalues. Rado's Theorem was introduced and
applied by Perfect in \cite{Perfect2}, to derive an important realizability
criterion for the \textit{RNIEP}. Surprisingly, this result was also ignored
in the literature about the problem until in \cite{Soto4}, the authors
rescue it and extend it to a new realizability criterion. Theorem \ref%
{TeoBrauerD} in Section $2,$ is a symmetric version of Theorem \ref%
{TeoBrauerC} and it was obtained in \cite{Soto6}. There, by the use of
Theorem \ref{TeoBrauerD} the authors give a criterion for the symmetric
realizability of a list $\Lambda =\{\lambda _{1},\ldots ,\lambda _{n}\}$ of
real numbers. This criterion, by its own definition, trivially contains any
other sufficient condition for the \textit{SNIEP} to have a solution.\\ \ \\
There are a number of known realizability criteria, which have
been obtained from the results of Brauer and/or Rado. Obviously they are not
included in this paper (see \cite{Borobia3, Perfect1, Perfect2, Soto, Soto4,
Soto7} for the \textit{RNIEP, }see\textit{\ }\cite{Soto3, Soto5, Soto6,
Soto7}\textit{\ }for the\textit{\ SNIEP, }see \cite{Borobia3, Rojo, Soto68,
Soto69} and the references therein for the complex case, and see \cite%
{Collao, Diaz, Julio, Soto65, Soto8, Soto9, Soto10, Soto11} for the \textit{%
URP}).
\\ \ \\ The paper is organized as follows: In Section 2 we introduce the
Theorems of Brauer and Rado above mentioned. In Section $3,$ from Brauer's
Theorem point of view, we prove some Guo's results \cite{Guo2}. In Section $4
$ we give, by applying Theorem \ref{TeoBrauerAA}, alternative proofs of
realizability criteria of Suleimanova \cite{Suleimanova} , Salzmann \cite%
{Salzmann}, Kellogg \cite{Kellogg}, Ciarlet \cite{Ciarlet}, Borobia \cite%
{Borobia1}, and \v{S}migoc \cite{Smigoc}. In Section $5$ we consider results
related to \textit{SNIEP. }In particular, from Theorem \ref{TeoBrauerD}, we
give a proof of two results due to Fiedler \cite{Fiedler}. In Section $6,$
we consider results associated with spectra of complex numbers, and we give,
from Rado's result point of view, a proof of a result due to \v{S}migoc \cite%
{Smigoc}. Finally, In Section $7,$ we examine the universal realizability
problem (\textit{URP}), and give an alternative proof of a Minc's result in
\cite{Minc3}.

\section{Brauer and Rado Theorems}

Brauer and Rado results, Theorem \ref{TeoBrauerAA} and Theorem \ref%
{TeoBrauerC}, respectively, have proven to be relevant for the study of the
\textit{NIEP} and the \textit{URP}. They have been applied with success to
generate sufficient conditions for the \textit{NIEP} and the \textit{URP} to
have a solution. These two theorems will be the unique results used along
this paper to give alternative proofs of distinct realizability criteria
compared in the maps constructed in \cite{Marijuan, Marijuan2}. We show that
virtually all known realizability criteria for the \textit{NIEP} and the
\textit{URP} can be obtained by applying Brauer or Rado results.

\begin{theorem}
\textit{[Brauer \cite{Brauer}]\label{TeoBrauerAA} Let }$A$\textit{\ be an }$%
n\times n$\textit{\ arbitrary matrix with eigenvalues }$\lambda _{1},\lambda
_{2},\ldots ,\lambda _{n}.$\textit{\ Let }$\mathbf{v}=(v_{1},v_{2},\ldots
,v_{n})^{T}$\textit{\ an eigenvector of }$A$\textit{\ associated with the
eigenvalue }$\lambda _{k}$\textit{\ and let }$\mathbf{q}$\textit{\ be any }$%
n $\textit{\ -dimensional vector. Then the matrix }$A+\mathbf{vq}^{T}$%
\textit{\ has eigenvalues }$\lambda _{1},\lambda _{2},\ldots ,\lambda
_{k-1},\lambda _{k}+\mathbf{v}^{T}\mathbf{q},\lambda _{k+1},\ldots ,\lambda
_{n}.$\textrm{\ }
\end{theorem}

\noindent Another proof, simpler that the one given in \cite{Brauer}, can be
found in \cite{Reams}. An immediate consequence of Theorem \ref{TeoBrauerAA}
is:

\begin{corollary}
\label{Epsilon}If $\Lambda =\{ \lambda _{1},\ldots ,\lambda_{n}\}$ is
realizable, then $\Lambda _{\epsilon }=\{ \lambda _{1}+\epsilon ,\lambda
_{2},\ldots ,\lambda _{n}\},$ $\epsilon \geq 0,$ is also realizable.
\end{corollary}

\begin{proof}
There exists a nonnegative matrix $A$ with spectrum $\Lambda ,$ which can be
taken as $A\in \mathcal{CS}_{\lambda _{1}}.$ Then, the matrix $A_{\epsilon
}=A+\epsilon \mathbf{ee}_{1}^{T}$ is nonnegative and, from Theorem \ref%
{TeoBrauerAA}, it has spectrum $\Lambda _{\epsilon }.$
\end{proof}

\ \newline
In \cite{Borobia3}, the authors introduce the concept of \textit{Brauer
negativity, }a quantity reflecting in a certain particular way how far $%
\Lambda $ is from being realized as the spectrum of a nonnegative matrix.
This negativity can be diminished by joining the list with a realizable
list, at best until the negativity is fully compensated and the joint list
becomes realizable. Then we have:

\begin{definition}
Given a list $\Lambda =\{\lambda _{1},\lambda _{2},\ldots ,\lambda _{n}\}$,
the \textit{Brauer negativity} of $\Lambda $ is
\begin{equation*}
\mathcal{N}(\Lambda )\equiv \min \{\delta \geq 0\,:\,\{\lambda _{1}+\delta
,\lambda _{2},\ldots ,\lambda _{n}\}\ \text{ is realizable}\}.\
\end{equation*}
\end{definition}
Note that a list $\Lambda $ is realizable if and only if $%
N(\Lambda )=0$.\\ \\ The following result is an extension of Theorem \ref%
{TeoBrauerAA}, and shows how to modify $r$ eigenvalues of a matrix of order $%
n,$ $r<n,$ via a $rank-r$ perturbation, without changing any of the $n-r$
remaining eigenvalues. This result was introduced by Perfect in \cite%
{Perfect2}. There, she point out that the result and its proof are due to R.
Rado.

\begin{theorem}
\textit{[Rado \cite{Perfect2}]\label{TeoBrauerC} Let A be an }$n\times n$%
\textit{\ arbitrary matrix with eigenvalues }$\lambda _{1},\lambda
_{2},\ldots ,\lambda _{n}.$\textit{\ Let }$X=[\mathbf{x}_{1}\mid \mathbf{x}%
_{2}\mid \cdots \mid \mathbf{x}_{r}]$\textit{\ be such that }$rank(X)=r$%
\textit{\ and }$A\mathbf{x}_{i}=\lambda _{i}\mathbf{x}_{i},$\textit{\ }$%
i=1,2,\ldots ,r,$\textit{\ }$r\leq n$\textit{. Let }$C$\textit{\ be an }$%
r\times n$\textit{\ arbitrary matrix. Then the matrix }$A+XC$\textit{\ has
eigenvalues }$\mu _{1},\mu _{2},\ldots ,\mu _{r},\lambda _{r+1},\lambda
_{r+2},\ldots ,\lambda _{n}$\textit{, where }$\mu _{1},\mu _{2},\ldots ,\mu
_{r}$\textit{\ are eigenvalues of the matrix}$\ \ \Omega +CX$\textit{\ with }%
$\Omega =diag\{\lambda _{1},\ldots ,\lambda _{r}\}$\textrm{$.$ }
\end{theorem}

\noindent Observe that for $r=1,$ Rado's Theorem is Brauer's Theorem. The
following example shows the importance of the Rado's result. In \cite{Soto4}%
, based on Theorem \ref{TeoBrauerC}, the authors give a sufficient
condition, which allow us, not only to decide on the realizability of the
list $\Lambda =\{6,3,3,-5,-5\},$ but also to construct a realizing matrix.
No one of known realizability criteria, other that the one given in \cite%
{Soto4}, is satisfied by $\Lambda $.

\begin{example}
\textit{\cite{Soto4} Consider }$\Lambda =\{6,3,3,-5,-5\}.$\textit{\ We
define the partition }$\Lambda _{0}=\{6,3,3\},$\textit{\ }$\Lambda
_{1}=\Lambda _{2}=\{-5\},$\textit{\ }$\Lambda _{3}=\emptyset ,$\textit{\
with the associated realizable lists }$\Gamma _{1}=\Gamma _{2}=\{5,-5\},$%
\textit{\ }$\Gamma _{3}=\{2\}.$\textit{\ The matrices }%
\begin{equation*}
A_{1}=A_{2}=%
\begin{bmatrix}
0 & 5 \\
5 & 0%
\end{bmatrix}%
,\ \ A_{3}=%
\begin{bmatrix}
2%
\end{bmatrix}%
\end{equation*}%
\textit{realize the lists }$\Gamma _{1}=\Gamma _{2}$\textit{\ and }$\Gamma
_{3},$\textit{\ respectively and the matrix }%
\begin{equation*}
A=%
\begin{bmatrix}
A_{1} & 0 & 0 \\
0 & A_{2} & 0 \\
0 & 0 & A_{3}%
\end{bmatrix}%
\end{equation*}%
\textit{has the spectrum }$\Gamma _{1}\cup \Gamma _{2}\cup \Gamma _{3}.$%
\textit{\ To apply the Rado's result we need to compute a }$3\times 3$%
\textit{\ nonnegative matrix with spectrum }$\Lambda _{0}$\textit{\ and
diagonal entries }$5,5,2.$\textit{\ From a Perfect's result \cite[Theorem $4$%
]{Perfect2} it is }%
\begin{equation*}
B=%
\begin{bmatrix}
5 & 0 & 1 \\
1 & 5 & 0 \\
0 & 4 & 2%
\end{bmatrix}%
\end{equation*}%
\textit{Then, for }%
\begin{equation*}
X=%
\begin{bmatrix}
1 & 0 & 0 \\
1 & 0 & 0 \\
0 & 1 & 0 \\
0 & 1 & 0 \\
0 & 0 & 1%
\end{bmatrix}%
,\ \ C=%
\begin{bmatrix}
0 & 0 & 0 & 0 & 1 \\
1 & 0 & 0 & 0 & 0 \\
0 & 0 & 4 & 0 & 0%
\end{bmatrix}%
,
\end{equation*}%
\textit{where the columns of }$X$\textit{\ are eigenvectors of }$A$\textit{\
and }$C$\textit{\ is obtained from }$B$\textit{\ in a certain appropriate
way (see \cite{Perfect2, Soto4}) we have that }%
\begin{equation*}
M=A+XC=%
\begin{bmatrix}
0 & 5 & 0 & 0 & 1 \\
5 & 0 & 0 & 0 & 1 \\
1 & 0 & 0 & 5 & 0 \\
1 & 0 & 5 & 0 & 0 \\
0 & 0 & 4 & 0 & 2%
\end{bmatrix}%
\end{equation*}%
\textit{is nonnegative with spectrum }$\Lambda $\textit{. }
\end{example}

\noindent In \cite{Soto6} the authors prove the following symmetric version
of Theorem \ref{TeoBrauerC}.

\begin{theorem}
\textit{\label{TeoBrauerD} \cite{Soto6} Let }$A$\textit{\ be an }$n\times n$%
\textit{\ symmetric matrix with eigenvalues }$\lambda _{1},\lambda
_{2},\ldots ,\lambda _{n}.$\textit{\ Let }$\{\mathbf{x}_{1},\mathbf{x}%
_{2},\ldots ,\mathbf{x}_{r}\}$\textit{\ an orthonormal set of eigenvectors
of }$A$\textit{\ such that }$AX=X\Omega ,$\textit{\ where }$X=[\mathbf{x}%
_{1}\mid \mathbf{x}_{2}\mid \ldots \mid \mathbf{x}_{r}]$\textit{\ and }$%
\Omega =diag\{\lambda _{1},\lambda _{2},\ldots ,\lambda _{r}\}.$\textit{\
Let }$C$\textit{\ be any }$r\times r$\textit{\ symmetric matrix. Then the
symmetric matrix }$A+XCX^{T}$\textit{\ has eigenvalues }$\mu _{1},\ldots
,\mu _{r},\lambda _{r+1},\ldots ,\lambda _{n}$\textit{, where }$\mu
_{1},\ldots ,\mu _{r}$\textit{\ are eigenvalues of the matrix }$\Omega +C.$%
\textrm{\ }
\end{theorem}

\section{On Guo results}

\noindent An important result of Guo \cite[Theorem 3.1]{Guo2}, which we
prove by using Theorem \ref{TeoBrauerC} establishes that:

\begin{theorem}
Let $\Lambda =\{ \lambda _{1},\lambda _{2},\ldots ,\lambda _{n}\},$ $\lambda
_{2}\in \mathbb{R},$ be realizable. Then for any $\epsilon >0,$ $\Lambda
_{\epsilon }=\{ \lambda _{1}+\epsilon ,\lambda _{2}\pm \epsilon ,\ldots
,\lambda _{n}\}$ is also realizable.
\end{theorem}

\begin{proof}
First we consider the case $\Lambda _{+\epsilon }=\{ \lambda _{1}+\epsilon
,\lambda _{2}+\epsilon ,\ldots ,\lambda _{n}\}.$ Let $\Lambda =\{ \lambda
_{1},\lambda _{2},\ldots ,\lambda _{n}\}$ be realizable with realizing
matrix $A\in \mathcal{CS}_{\lambda _{1}}.$ Then $A\mathbf{e}=\lambda _{1}%
\mathbf{e.}$ Let $A\mathbf{x}=\lambda _{2}\mathbf{x,}$ with $\mathbf{x}%
^{T}=\left( x_{1},x_{2},\ldots ,x_{n}\right) ,$ $x_{1}=\max
\{x_{1},x_{2},\ldots ,x_{n}\},$ $x_{2}=\min \{x_{1},x_{2},\ldots ,x_{n}\}.$
It is clear that $x_{1}\geq 0$ and $x_{2}\leq 0.$ Let

\begin{equation*}
X=%
\begin{bmatrix}
1 & x_{1} \\
1 & x_{2} \\
\vdots & \vdots \\
1 & x_{n}%
\end{bmatrix}%
,\ C=%
\begin{bmatrix}
c_{11} & c_{12} & \cdots & \cdots & c_{1n} \\
c_{21} & c_{22} & \cdots & \cdots & c_{2n}%
\end{bmatrix}%
,\text{ and }\Omega =%
\begin{bmatrix}
\lambda _{1} & 0 \\
0 & \lambda _{2}%
\end{bmatrix}%
,
\end{equation*}%
with
\begin{eqnarray*}
c_{11} &=&\frac{-\epsilon x_{2}}{x_{1}-x_{2}},\ \ c_{12}=\frac{\epsilon x_{1}%
}{x_{1}-x_{2}},\ \ c_{1j}=0,\text{ \ }j=3,\ldots ,n \\
&& \\
c_{21} &=&\frac{\epsilon }{x_{1}-x_{2}},\ \ c_{22}=\frac{-\epsilon }{%
x_{1}-x_{2}},\text{ \ }c_{2j}=0,\text{ \ }j=3,\ldots ,n.
\end{eqnarray*}%
Then
\begin{equation*}
XC=%
\begin{bmatrix}
c_{11}+c_{21}x_{1} & c_{12}+c_{22}x_{1} & 0 & \cdots & 0 \\
c_{11}+c_{21}x_{2} & c_{12}+c_{22}x_{2} & 0 & \cdots & 0 \\
\vdots & \vdots & \vdots & \ddots & \vdots \\
c_{11}+c_{21}x_{n} & c_{12}+c_{22}x_{n} & 0 & \cdots & 0%
\end{bmatrix}%
\geq 0\text{ and }CX=%
\begin{bmatrix}
\epsilon & 0 \\
0 & \epsilon%
\end{bmatrix}%
.
\end{equation*}%
Now, from Rado's Theorem, $A+XC$ is nonnegative with spectrum $\Lambda
_{+\epsilon }=\{ \lambda _{1}+\epsilon ,\lambda _{2}+\epsilon ,\ldots
,\lambda _{n}\}.$ For the case $\Lambda _{-\epsilon }=\{ \lambda
_{1}+\epsilon ,\lambda _{2}-\epsilon ,\ldots ,\lambda _{n}\}$ we take
\begin{eqnarray*}
c_{11} &=&\frac{\epsilon x_{1}}{x_{1}-x_{2}},\ \ c_{12}=\frac{-\epsilon x_{2}%
}{x_{1}-x_{2}},\ \ c_{1j}=0,\text{ \ }j=3,\ldots ,n \\
&& \\
c_{21} &=&\frac{-\epsilon }{x_{1}-x_{2}},\ \ c_{22}=\frac{\epsilon }{%
x_{1}-x_{2}},\text{ \ }c_{2j}=0,\text{ \ }j=3,\ldots ,n.
\end{eqnarray*}%
and the result follows.
\end{proof}

\ \newline
A weaker version of this result is:

\begin{lemma}
Let $\Lambda _{1}=\{ \alpha _{1},\alpha _{2},\ldots ,\alpha _{n}\}$ and $%
\Lambda _{2}=\{ \beta _{1},\beta _{2},\ldots ,\beta _{m}\}$ be realizable
lists. Then for any $\epsilon \geq \max \{ \beta _{1}-\alpha _{1},0\},$
\begin{equation*}
\Lambda =\{ \alpha _{1}+\epsilon ,\beta _{1}-\epsilon ,\alpha _{2},\ldots
,\alpha _{n},\beta _{2},\ldots ,\beta _{m}\} \text{ is realizable.}
\end{equation*}
\end{lemma}

\begin{proof}
There exist $A_{1}\in \mathcal{CS}_{\alpha _{1}}$ realizing $\Lambda _{1}$
and $A_{2}\in \mathcal{CS}_{\beta _{1}}$ realizing $\Lambda _{2}.$ Then
\begin{equation*}
(A_{2}-\epsilon \mathbf{ee}_{1}^{T})\in \mathcal{CS}_{\beta _{1}-\epsilon }
\end{equation*}
and
\begin{equation*}
B=
\begin{bmatrix}
A_{2}-\epsilon \mathbf{ee}_{1}^{T} & (\alpha _{1}-\beta _{1}+\epsilon )%
\mathbf{ee}_{1}^{T} \\
0 & A_{1}%
\end{bmatrix}
\in \mathcal{CS}_{\alpha _{1}}
\end{equation*}
Then $A=(B+\epsilon \mathbf{ee}_{1}^{T})\in \mathcal{CS}_{\alpha
_{1}+\epsilon }$ is nonnegative with spectrum $\Lambda .$
\end{proof}

\bigskip

\noindent A result due to Guo \cite[Theorem $2.1$]{Guo2}, states the
existence of a real number $\lambda _{0},$

\begin{equation}
\max_{2\leq j\leq n}\left \vert \lambda _{j}\right \vert \leq \lambda
_{0}\leq 2n\max_{2\leq j\leq n}\left \vert \lambda _{j}\right \vert ,
\label{ec1}
\end{equation}
such that the list of complex numbers $\Lambda =\{ \lambda _{1},\lambda
_{2},\ldots,\lambda _{n}\}$ is realizable if and only if $\lambda _{1}\geq
\lambda _{0}.$ The Guo index $\lambda _{0}$ is the minimum $\lambda $ such
that $\{ \lambda ,\lambda _{2},\ldots,\lambda _{n}\}$ is realizable. The
problem of finding $\lambda _{0}$ is not solved in the paper of Guo. Here,
we show that the upper bound in (\ref{ec1}) may be reduced to $%
(n-1)\max_{2\leq j\leq n}\left \vert \lambda _{j}\right \vert ,$ and that
this bound is sharp.

\begin{theorem}
Let $\Lambda ^{\prime }=\{ \lambda _{2},\ldots,\lambda _{n}\} \subset
\mathbb{C} $ be such that $\overline{\Lambda ^{\prime }}=\Lambda ^{\prime }$
with $\lambda _{2}\geq \cdots \geq \lambda _{p}$ real, $2\leq p\leq n.$
Then, the Guo upper bound in (\ref{ec1}) may be reduced to$\ $
\begin{equation}
\lambda _{0}\leq (n-1)\max_{2\leq j\leq n}\left \vert \lambda _{j}\right
\vert .  \label{ec2.2}
\end{equation}
\end{theorem}

\begin{proof}
Let $m=\max_{2\leq j\leq n}\left \vert \lambda _{j}\right \vert $ and let $%
\mu _{j}=\frac{\lambda _{j}}{m(n-1)},$ $j=2,3,\ldots,n.$ Then, $\Gamma
^{\prime }=\{ \mu _{2},\ldots ,\mu _{n}\}$ is a list of complex numbers such
that $\ \left \vert \mu _{j}\right \vert \leq \frac{1}{n-1},$ $j=2,3,\ldots
,n.$ Consider the initial matrix
\begin{equation*}
B=%
\begin{bmatrix}
0 & 0 & 0 & \cdots & \cdots & \cdots & \cdots & \cdots & 0 \\
-\mu _{2} & \mu _{2} & \ddots &  &  &  &  &  & \vdots \\
\vdots & \ddots & \ddots & \ddots &  &  &  &  & \vdots \\
-\mu _{p} & \vdots & \ddots & \mu _{p} & 0 &  &  &  & \vdots \\
-x_{s} & y_{s} &  & \ddots & x_{s} & -y_{s} &  &  & \vdots \\
-x_{s} & -y_{s} &  &  & y_{s} & x_{s} & \ddots &  & \vdots \\
\vdots & \vdots &  &  &  & \ddots & \ddots & \ddots & 0 \\
-x_{t} & y_{t} &  &  &  &  & \ddots & x_{t} & -y_{t} \\
-x_{t} & -y_{t} & \cdots & \cdots & \cdots & \cdots & 0 & y_{t} & x_{t}%
\end{bmatrix}%
,
\end{equation*}
where $\mu _{2},\mu _{3},\ldots ,\mu _{p}$ are real, $x_{j}=Re\mu _{j},$ $%
y_{j}=Im\mu _{j},$ $p+1\leq j\leq \frac{n+p}{2}.$ Then $B\in \mathcal{CS}%
_{0} $ has eigenvalues $0,\mu _{2},\ldots ,\mu _{p},\mu _{p+1},\ldots ,\mu
_{n}.$ \newline
If $Re\mu _{j}\leq 0,$ $j=2,3,\ldots ,n,$ then all the entries in the first
column of $B$ are nonnegative. Let
\begin{equation*}
\mathbf{q}=(0,\frac{1}{n-1},\ldots ,\frac{1}{n-1})^{T}.
\end{equation*}
From Theorem \ref{TeoBrauerAA}, the matrix $A^{\prime }=B+\mathbf{eq}^{T}$
is nonnegative with eigenvalues $1,\mu _{2},\ldots ,\mu _{n}$ and the matrix
$A=m(n-1)A^{\prime }$ is nonnegative with eigenvalues $(n-1)m,\lambda
_{2},\ldots ,\lambda _{n}.$\newline
If $Re\mu _{k}>0$ for some $k,$ $3\leq k\leq n$, then all the entries in the
$k^{th}$ column of $B$ (or in the $(k-1)^{th}$ column of $B$ if $k$
corresponds to the second column in the corresponding $2\times 2$ complex
block) are nonnegative. Let
\begin{equation*}
\mathbf{q}=(\frac{1}{n-1},\ldots ,\frac{1}{n-1},0,\frac{1}{n-1},\ldots,
\frac{1}{n-1})^{T}
\end{equation*}
with the zero in the $k^{th}$ position ($(k-1)^{th}$ position). Then the
matrix $A^{\prime }=B+\mathbf{eq}^{T}$ is nonnegative with eigenvalues $%
1,\mu _{2},\mu _{3,}\ldots ,\mu _{n}$ and the matrix $A=m(n-1)A^{\prime }$
is nonnegative with eigenvalues $(n-1)m,\lambda _{2},\lambda _{3},\ldots
,\lambda _{n}.$ \newline
If $\mu _{2}>0$ with $Re\mu _{j}<0,$ $j=3,\ldots ,n,$ then we write the $%
-Re\mu _{j}{\acute{}}s,$ \ $3\leq j\leq n,$ along the second column of $B$
and the $\pm Im\mu _{j}{\acute{}}s,$ $p+1\leq j\leq n,$ along the first
column of $B.$ Then, with
\begin{equation*}
\mathbf{q}=(\frac{1}{n-1},0,\frac{1}{n-1},\ldots ,\frac{1}{n-1})^{T}
\end{equation*}
we obtain as before, the nonnegative matrix $A=m(n-1)A^{\prime }$ with the
required eigenvalues.
\end{proof}

\bigskip

\noindent The inequality (\ref{ec2.2}) is sharp. In fact, let $\lambda
_{j}=-1,$ $2\leq j\leq n.$ Then, from a Suleimanova result \cite{Suleimanova}
(Theorem \ref{Sule} in this paper), the problem has a solution if and only
if \ $\lambda _{1}\geq \lambda _{0}=(n-1).$ Thus, in the real case, the Guo
result guarantees the existence of a nonnegative matrix $A$ with spectrum $%
\Lambda =\{\lambda _{1},\lambda _{2},\ldots ,\lambda _{n}\}$ for all
\begin{equation}
\lambda _{1}\geq (n-1)\max_{2\leq j\leq n}\left\vert \lambda _{j}\right\vert
.  \label{ec3}
\end{equation}%
If $\Lambda =\{\lambda _{1},\lambda _{2},\ldots ,\lambda _{n}\},$ with $%
\sum\limits_{i=1}^{n}\lambda _{i}=0,$ is a realizable list of complex
numbers, then $\lambda _{1}$ is the Guo index.

\begin{remark}
\textit{In \cite{Rojo} the authors show how to compute the Guo index }$%
\lambda _{0}$\textit{\ for circulant nonnegative matrices. There, they give
a necessary and sufficient condition for a list of complex numbers }$\Lambda
=\{\lambda _{1},\lambda _{2},\ldots ,\lambda _{n}\}$\textit{\ to be the
spectrum of a circulant nonnegative matrix. However, to compute }$\lambda
_{0}$\textit{\ becomes a prohibitive task for large }$n.$\textit{\ Then,
they prove, by the use of Theorem \ref{TeoBrauerAA}, a more handleable
realizability criterion (see \cite{Rojo}).}\textrm{\ }
\end{remark}

\section{RNIEP}

\noindent In this section we consider the real nonnegative inverse
eigenvalue problem, and we show that the realizability criteria given by
Suleimanova \cite{Suleimanova}, Salzmann \cite{Salzmann}, Ciarlet \cite%
{Ciarlet}, Kellogg \cite{Kellogg}, and Borobia \cite{Borobia1}, may all be
obtained by applying the Brauer's result. Perfect \cite{Perfect1} was the
first one who used Theorem \ref{TeoBrauerAA} to derive sufficient conditions
for the \emph{RNIEP} to have a solution. There are a number of realizability
criteria, which are not considered here because they were obtained from
Brauer's result or Rado's result.\

\bigskip

\noindent We also examine, in this section, a result due to \v{S}migoc \cite[%
Theorem $10$]{Smigoc}, which can be proved by using Theorem \ref{TeoBrauerC}
. Since Rado's result is a generalization of Brauer's result, criteria from
Rado give, in general, better information about the realizability of a given
real list $\Lambda =\{\lambda _{1},\ldots ,\lambda _{n}\}.$ All proofs of
criteria generate from Brauer's result or from Rado's result are
constructive, in the sense that they always allow us to compute a realizing
matrix.

\begin{theorem}
\textit{[Suleimanova \cite{Suleimanova}]\label{Sule} Let }$\Lambda
=\{\lambda _{1},\lambda _{2},\cdots ,\lambda _{n}\}$\textit{\ be a list of
real numbers satisfying }$\lambda _{1}+\lambda _{2}+\ldots +\lambda _{n}\geq
0,\text{ \ }\lambda _{k}<0,\ k=2,3,\ldots ,n.$\textit{\ Then }$\Lambda $%
\textit{\ is realizable. }
\end{theorem}

\begin{proof}
It was proved in \cite{Soto} by the use of Brauer's Theorem \ref{TeoBrauerAA}%
.
\end{proof}

\begin{theorem}
\textit{[Salzmann \cite{Salzmann}]\label{TeoSalzmann} Let }$\Lambda
=\{\lambda _{1},\lambda _{2},\ldots ,\lambda _{n}\}$\textit{\ be such that }%
\begin{equation}
\sum_{k=1}^{n}\lambda _{k}\geq 0,  \label{EcSalz1}
\end{equation}%
\textit{and }%
\begin{equation}
\lambda _{k}+\lambda _{n-k+1}\leq \frac{2}{n}\sum_{k=1}^{n}\lambda _{k},%
\text{ \ }k=2,3,\ldots ,\left[ \frac{n+1}{2}\right] .  \label{EcSalz2}
\end{equation}%
\textit{Then }$\Lambda $\textit{\ is realizable. }
\end{theorem}

\begin{proof}
We only need to prove the assertion for $\Lambda $ satisfying $%
\sum_{k=1}^{n}\lambda _{k}=0$. Therefore, suppose conditions (\ref{EcSalz1})
and (\ref{EcSalz2}) are satisfied with $\sum_{k=1}^{n}\lambda _{k}=0.$ Then $%
S_{k}=\lambda _{k}+\lambda _{n-k+1}\leq 0$ for $k=2,3,\ldots ,\left[ \frac{n%
}{2}\right] $\ and $S_{\frac{n+1}{2}}=\lambda _{\frac{n+1}{2}}\leq 0$\ for $%
n $ odd. Now, we apply \cite[Theorem $11$]{Soto}, which is a criterion
obtained by applying Theorem \ref{TeoBrauerAA}. Let
\begin{equation*}
B=%
\begin{bmatrix}
B_{11} & 0 & \ddots & 0 \\
B_{21} & B_{22} & \ddots & \ddots \\
\ddots & \ddots & \ddots & 0 \\
B_{\left[ \frac{n}{2}\right] 1} & \ddots & 0 & B_{\left[ \frac{n}{2}\right] }%
\end{bmatrix}%
\end{equation*}
with diagonal blocks
\begin{equation*}
B_{11}=%
\begin{bmatrix}
0 & 0 \\
-\lambda _{n} & \lambda _{n}%
\end{bmatrix}
,\text{ }B_{kk}=
\begin{bmatrix}
0 & \lambda _{k} \\
-\lambda _{n-k+1} & \lambda _{k}+\lambda _{n-k+1}%
\end{bmatrix}%
\end{equation*}
$\ $ and
\begin{equation*}
B_{k1}=
\begin{bmatrix}
0 & -\lambda _{k} \\
0 & -\lambda _{k}%
\end{bmatrix}
,\text{ }k=2,3,\ldots ,\left[ \frac{n}{2}\right]
\end{equation*}
Let $\mathbf{q}=(q_{1},q_{2},\ldots ,q_{n})^{T}$ \ with
\begin{eqnarray*}
q_{2k-1} &=&0,\text{ \ }k=1,2,\ldots ,\frac{n}{2}, \\
q_{2} &=&-\lambda _{n}, \\
q_{2k} &=&-(\lambda _{k}+\lambda _{n-k+1}),\text{ \ }k=2,3,\ldots ,\frac{n}{2%
}.
\end{eqnarray*}%
Then, $A=B+\mathbf{eq^{T}}$ is nonnegative with spectrum $\Lambda .$
\end{proof}

\begin{theorem}
\textit{[Ciarlet \cite{Ciarlet}] Let \ }$\Lambda =\{\lambda _{1},\lambda
_{2},\ldots ,\lambda _{n}\}$\textit{\ such that }%
\begin{equation*}
\lambda _{1}\geq \lambda _{2}\geq \cdots \geq \lambda _{p}\geq 0>\lambda
_{p+1}\geq \cdots \geq \lambda _{n}\ \ \ \ \ \ 1\leq p\leq n-1.
\end{equation*}%
\textit{If }$|\lambda _{k}|\leq \frac{\lambda _{1}}{n},\ k=2,3,\ldots ,n$%
\textit{, then }$\Lambda $\textit{\ is realizable. }
\end{theorem}

\begin{proof}
Consider the matrix%
\begin{equation*}
B=
\begin{bmatrix}
0 & 0 & 0 & \cdots & 0 & 0 & \cdots & 0 \\
-\lambda _{2} & \lambda _{2} & 0 & \cdots & 0 & 0 & \cdots & 0 \\
-\lambda _{3} & 0 & \lambda _{3} & \cdots & 0 & 0 & \cdots & 0 \\
\vdots & \vdots & \vdots & \ddots & \vdots & \vdots & \ddots & \vdots \\
-\lambda _{p} & 0 & 0 & \cdots & \lambda _{p} & 0 & \cdots & 0 \\
-\lambda _{p+1} & 0 & 0 & \cdots & 0 & \lambda _{p+1} & \cdots & 0 \\
\vdots & \vdots & \vdots & \cdots & \vdots & \vdots & \ddots & \vdots \\
-\lambda _{n} & 0 & 0 & \cdots & 0 & 0 & \cdots & \lambda _{n}%
\end{bmatrix}
\in \mathcal{CS}_{0}.
\end{equation*}
$B$ has spectrum $\{0,\lambda _{2},\ldots ,\lambda _{n}\}$. Then, for $%
\mathbf{q}=(\frac{\lambda _{1}}{n},\ldots ,\frac{\lambda _{1}}{n})^{T},$ we
have that $A=B+\mathbf{eq^{T}}$ has spectrum $\Lambda ,$ and since $|\lambda
_{k}|\leq \frac{\lambda _{1}}{n},\ \ k=2,3,\ldots ,n$, $A$ is nonnegative.
\end{proof}

\begin{theorem}
\textit{[Kellogg \cite{Kellogg}]\label{TeoKe} Let }$\Lambda =\{\lambda
_{1},\lambda _{2},\ldots ,\lambda _{n}\}$\textit{\ be a list of real numbers
with }$\lambda _{1}\geq \lambda _{2}\geq \ldots \geq \lambda _{n}$\textit{\
and let }$p$\textit{\ be the greatest index }$j$\textit{\ }$(1\leq j\leq n)$%
\textit{\ for which }$\lambda _{j}\geq 0.$\textit{\ Let the set of indices }%
\begin{equation*}
K=\{i:\lambda _{i}\geq 0\text{ and }\lambda _{i}+\lambda _{n-i+2}<0,\text{ \
}i\in \{2,3,\ldots ,\left[ \frac{n+1}{2}\right] \}\}.
\end{equation*}%
\textit{If }%
\begin{equation}
\lambda _{1}\geq -\sum_{i\in K,\text{ }i<k}(\lambda _{i}+\lambda
_{n-i+2})-\lambda _{n-k+2}\text{\ \ \ for all }k\in K,  \label{Kec1}
\end{equation}%
$\ $\textit{and }%
\begin{equation}
\lambda _{1}\geq -\sum_{i\in K}(\lambda _{i}+\lambda
_{n-i+2})-\sum_{j=p+1}^{n-p+1}\lambda _{j},\ \ \text{provided that}\ n\geq
2p,  \label{Kec2}
\end{equation}%
\textit{then }$\Lambda $\textit{\ is realizable. }
\end{theorem}

\begin{proof}
Suppose conditions (\ref{Kec1}) and ( \ref{Kec2}) are satisfied. Let $%
K=\{k_{1},k_{2},\ldots ,k_{t}\}$ be the Kellogg set of indices. Consider the
partition $\Lambda =\Lambda _{1}\cup \cup _{i=1}^{t}\Lambda _{k_{i}}\cup
\Lambda _{\mathcal{R}},$ where
\begin{eqnarray*}
\Lambda _{1} &=&\{\lambda _{1},\lambda _{p+1},\ldots ,\lambda _{n-p+1}\} \\
\Lambda _{k_{i}} &=&\{\lambda _{k_{i}},\lambda _{n-k_{i}+2}\},\text{ }%
k_{i}\in K\ i=1,2,\ldots ,t\ \ \ (\lambda _{k_{1}}\geq \lambda _{k_{2}}\geq
\cdots \geq \lambda _{k_{t}}\geq 0) \\
\Lambda _{\mathcal{R}} &=&\Lambda -\Lambda _{1}-\cup _{i=1}^{t}\Lambda
_{k_{i}}.
\end{eqnarray*}%
From \eqref{Kec2} we have that
\begin{equation*}
\lambda _{1}+\sum_{j=p+1}^{n-p+1}\lambda _{j}\geq -\sum_{\substack{ i\in K}}%
(\lambda _{i}+\lambda _{n-i+2})>0\ \ ((\lambda _{i}+\lambda _{n-i+2})<0\
\forall i\in K).
\end{equation*}%
Since $\lambda _{p+1},\lambda _{p+2},\ldots ,\lambda _{n-p+1}$ are negative,
then from Theorem \ref{Sule}, $\Lambda _{1}$ is realizable.\ \newline
\ \newline
If $\Lambda _{\mathcal{R}}=\emptyset $ then $\Lambda =\Lambda _{1}\ \cup \
\cup _{i=1}^{t}\Lambda _{k_{i}}.$ If $\Lambda _{\mathcal{R}}\neq \emptyset $%
, then $\Lambda _{\mathcal{R}}$ contain sublists $\Lambda _{i}=\{\lambda
_{i},\lambda _{n-i+2}\}$ such that $(\lambda _{i}+\lambda _{n-i+2})\geq 0$.
So, the lists $\Lambda _{i}$ are realizable by
\begin{equation*}
A_{i}=%
\begin{bmatrix}
0 & \lambda _{i} \\
-\lambda _{n-i+2} & \lambda _{i}+\lambda _{n-i+2}%
\end{bmatrix}%
\end{equation*}%
Let $A_{\mathcal{R}}$ be the realizing matrix of $\Lambda _{\mathcal{R}}$.
On the other hand, for each sublist $\Lambda _{k_{i}}=\{\lambda
_{k_{i}},\lambda _{n-k_{i}+2}\}$ with $k_{i}\in K$, $i=1,2,\ldots ,t,$
\begin{equation*}
B_{k_{i}}=%
\begin{bmatrix}
0 & \lambda _{k_{i}} \\
-\lambda _{n-k_{i}+2} & \lambda _{k_{i}}+\lambda _{n-k_{i}+2}%
\end{bmatrix}%
\end{equation*}%
is a $2\times 2$ matrix with spectrum $\Lambda _{k_{i}}$ \ $i=1,2,\ldots ,t$%
. Since $\lambda _{k_{i}}\geq 0$ and $(\lambda _{k_{i}}+\lambda
_{n-k_{i}+2})<0$, then
\begin{equation*}
-\lambda _{n-k_{i}+2}>0.
\end{equation*}%
Let
\begin{equation*}
\Lambda _{1}^{\prime }=\{\lambda _{1}+\sum_{i=1}^{t}(\lambda
_{k_{i}}+\lambda _{n-k_{i}+2}),\lambda _{p+1},\ldots ,\lambda _{n-p+1}\},
\end{equation*}%
From \eqref{Kec2} $\Lambda _{1}^{\prime }$ is realizable by an $%
(n-2p+2)\times (n-2p+2)$ matrix $B_{1}^{\prime }\in \mathcal{CS}_{\mu },$
where
\begin{equation*}
\mu =\lambda _{1}+\sum_{i=1}^{t}(\lambda _{k_{i}}+\lambda _{n-k_{i}+2}).
\end{equation*}%
Now, let
\begin{equation*}
B=%
\begin{bmatrix}
B_{1}^{\prime } & 0 & 0 & \cdots  & \cdots  & 0 & 0 \\
B_{2} & B_{k_{t}} & 0 & \ddots  & \vdots  & \vdots  & \vdots  \\
0 & B_{k_{t-1}}(\mu ) & B_{k_{t-1}} & \ddots  & \ddots  & \vdots  & \vdots
\\
\vdots  & B_{k_{t-2}}(\mu ) & B_{k_{t-1}}^{\prime } & B_{k_{t-2}} & \ddots
& \ddots  & \vdots  \\
\vdots  & \vdots  & \vdots  & \vdots  & \ddots  & 0 & 0 \\
\vdots  & B_{k_{2}}(\mu ) & B_{k_{t-1}}^{\prime } & B_{k_{t-2}}^{\prime } &
\cdots  & B_{k_{2}} & 0 \\
0 & B_{k_{1}}(\mu ) & B_{k_{t-1}}^{\prime } & B_{k_{t-2}}^{\prime } & \cdots
& B_{k_{2}}^{\prime } & B_{k_{1}}%
\end{bmatrix}%
,
\end{equation*}%
where
\begin{equation*}
B_{2}=%
\begin{bmatrix}
0 & \cdots  & 0 & \mu -\lambda _{k_{t}} \\
0 & \cdots  & 0 & \mu -\lambda _{k_{t}}%
\end{bmatrix}%
_{2\times (n-2p+2)},
\end{equation*}

\begin{equation*}
B_{k_{i}}=%
\begin{bmatrix}
0 & \lambda _{k_{i}} \\
-\lambda _{n-k_{i}+2} & \lambda _{k_{i}}+\lambda _{n-k_{i}+2}%
\end{bmatrix}
\ \ \ \ \ \ i=1,2,\ldots ,t,
\end{equation*}
\begin{equation*}
B_{k_{i}}(\mu )=
\begin{bmatrix}
0 &  & \mu -\lambda _{k_{i}}-\sum \limits_{j=i+1}^{t-1}(\lambda
_{k_{j}}+\lambda _{n-k_{j}+2}) \\
0 &  & \mu -\lambda _{k_{i}}-\sum \limits_{j=i+1}^{t-1}(\lambda
_{k_{j}}+\lambda _{n-k_{j}+2})%
\end{bmatrix}
\ \ \ \ \ i=1,2,\cdots ,t-1,
\end{equation*}
\begin{equation*}
B_{k_{i}}^{\prime }=
\begin{bmatrix}
0 & \lambda _{k_{i}}+\lambda _{n-k_{i}+2} \\
0 & \lambda _{k_{i}}+\lambda _{n-k_{i}+2}%
\end{bmatrix}
\ \ \ \ \ i=2,3,\ldots ,t-1
\end{equation*}
are such that $B\in \mathcal{CS}_{\mu }$ with spectrum $\Lambda _{1}^{\prime
}\cup \cup _{i=1}^{t}\Lambda _{k_{i}}$. \newline
Note that the sum $\sum \limits_{j=i+1}^{t-1}(\lambda _{k_{j}}+\lambda
_{n-k_{j}+2})$ \ in the matrix $B_{k_{i}}(\mu )$ is zero if $i+1>t-1$ for
some $i$. From Brauer Theorem \ref{TeoBrauerAA}, we have for
\begin{equation*}
\mathbf{q}=(\underset{(n-2p+2)-times}{\underbrace{0,\ldots ,0}},0,-(\lambda
_{k_{t}}+\lambda _{n-k_{t}+2}),\ldots ,0,-(\lambda _{k_{1}}+\lambda
_{n-k_{1}+2}))^{T}\ ,
\end{equation*}%
that $M=B+\mathbf{eq^{T}}\in \mathcal{CS}_{\lambda _{1}}$ is nonnegative
with spectrum $\Lambda _{1}\cup \cup _{i=1}^{t}\Lambda _{k_{i}}.$ \ \newline
Finally, a matrix $A=M\oplus A_{R}$ is nonnegative with spectrum $\Lambda $.
\end{proof}

\ \newline
Before to prove Borobia result by using Theorem \ref{TeoBrauerAA}, we need
the following result, which was proved in \cite{Soto2}:

\begin{lemma}
\label{SotoBoro}Let
\begin{equation*}
B_{k}=
\begin{bmatrix}
0 & \lambda _{k} \\
-\lambda _{n-k+2} & \lambda _{k}+\lambda _{n-k+2}%
\end{bmatrix}
,
\end{equation*}
where, $\lambda _{n-k+2}=\sum_{j=1}^{r}\mu _{j}<0$ with $\mu _{j}<0,$ $\
j=1,2,\ldots ,r.$ Then
\begin{equation*}
\Lambda _{k}=\{ \lambda _{k},\lambda _{n-k+2}\} \text{ \ and \ }\Lambda
_{r}=\{ \lambda _{k},\mu _{1},\ldots ,\mu _{r}\}
\end{equation*}%
have the same Brauer negativity, $\mathcal{N}(\Lambda _{k})=\mathcal{N}%
(\Lambda _{r}).$ Moreover, there exists a matrix $B_{r}\in \mathcal{CS}%
_{\lambda _{k}},$ of order $r+1,$ with spectrum $\Lambda _{r}.$
\end{lemma}

\begin{theorem}
\textit{[Borobia, \cite{Borobia1}]\label{TeoBoro} Let }$\Lambda =\{\lambda
_{1},\lambda _{2},\ldots ,\lambda _{n}\}$\textit{\ be a list of real numbers
with }$\lambda _{1}\geq \lambda _{2}\geq \ldots \geq \lambda _{p}\geq
0>\lambda _{p+1}\geq \cdots \geq \lambda _{n}$\textit{\ (}$p$\textit{\ is
the greatest index }$j$\textit{\ }$(1\leq j\leq n)$\textit{\ for which }$%
\lambda _{j}\geq 0.$\textit{) If there exists a partition }$J_{1}\cup
J_{2}\cup \ldots \cup J_{t}$\textit{\ of }$J=\{\lambda _{p+1},\lambda
_{p+2},\ldots ,\lambda _{n}\},$\textit{\ for some }$1\leq t\leq n-p,$\textit{%
\ such that }%
\begin{equation}
\lambda _{1}\geq \lambda _{2}\geq \ldots \geq \lambda _{p}\geq \sum_{\lambda
\in J_{1}}\lambda \geq \sum_{\lambda \in J_{2}}\lambda \geq \ldots \geq
\sum_{\lambda \in J_{t}}\lambda  \label{Boroec1}
\end{equation}%
\textit{satisfies the Kellogg conditions (\ref{Kec1}) and (\ref{Kec2}), then
}$\Lambda $\textit{\ is realizable. }
\end{theorem}

\begin{proof}
Suppose that Borobia realizability criterion is satisfied. Then, there
exists a partition
\begin{equation*}
J_{1}\cup J_{2}\cdots \cup J_{t}\ \text{ of}\ J=\{\lambda _{p+1},\ldots
,\lambda _{n}\},
\end{equation*}%
with
\begin{equation*}
\sum_{\lambda \in J_{j}}\lambda =\mu _{p+j},\ \ j=1,2,\ldots ,t,\ \ 1\leq
t\leq n-p,
\end{equation*}%
such that the new list
\begin{equation*}
\Gamma =\{\mu _{1},\mu _{2},\ldots ,\mu _{p},\mu _{p+1},\ldots ,\mu _{p+t}\}
\end{equation*}%
with $\mu _{1}\geq \cdots \geq \mu _{p}\geq 0>\mu _{p+1}\geq \cdots \geq \mu
_{p+t}$, where $\mu _{i}=\lambda _{i},$ \ $i=1,\ldots ,p,$ satisfies the
Kellogg realizability criterion.\ \newline
Now, we apply the same proof as in Theorem \ref{TeoKe} of Kellogg, except
for one detail: Since the new list
\begin{equation*}
\Gamma =\{\mu _{1},\ldots ,\mu _{p},\mu _{p+1},\ldots ,\mu _{p+t}\}
\end{equation*}%
has less elements than the original list $\Lambda =\{\lambda _{1},\lambda
_{2},\ldots ,\lambda _{n}\}$ and we want to obtain a nonnegative matrix of
order $n$ realizing $\Lambda ,$ then in each step of the proof, before to
manipulate any of the new sublists
\begin{equation*}
\Gamma _{ki}=\{\mu _{k_{i}},\mu _{p+t-k_{i}+2}\},\text{ }k_{i}\in K\ \text{%
(Kellogg set of indices)}
\end{equation*}%
where,
\begin{equation*}
\mu _{p+t-k_{i}+2}=\mu _{p+j}=\sum_{\lambda \in J_{j}}\lambda <0,\ \forall
k_{i}\in K\ \text{and}\ \ j=1,\ldots ,t.
\end{equation*}%
We must extend $\Gamma _{ki}$ to the lists $\Gamma _{r_{i}}=\{\mu
_{k_{i}},\lambda _{1_{j}},\ldots ,\lambda _{r_{j}}\},\forall k_{i}\in K,$
where
\begin{equation*}
\lambda _{s_{j}}\in J_{j},s=1,\ldots ,r,\ \text{and}\ \sum_{s=1}^{r}\lambda
_{s_{j}}=\mu _{p+t-k_{i}+2},\ j=1,\ldots ,t.
\end{equation*}%
From Lemma \ref{SotoBoro} $\mathcal{N}(\Gamma _{k_{i}})=\mathcal{N}(\Gamma
_{r_{i}}).$ Then the realizing matrix for $\Gamma _{r_{i}}$ is of the
required size. In the same way, if necessary, we must also extend the new
list
\begin{equation*}
\Gamma _{1}^{\prime }=\{\mu ,\mu _{p+1},\ldots ,\mu _{t+1}\},\ \text{with}\
\mu =\mu _{1}+\sum_{k_{i}\in K}(\mu _{k_{i}}+\mu _{p+t-k_{i}+2})
\end{equation*}%
to a new lits $\Gamma _{1}^{\prime \prime }$ with the same Brauer
negativity, by replacing the corresponding $\mu _{p+1},\ldots ,\mu _{t+1}$
by $\sum_{\lambda \in J_{1}}\lambda ,\ldots ,\sum_{\lambda \in
J_{t-p+1}}\lambda ,$ respectively. Again from Lemma \ref{SotoBoro} $\mathcal{%
N}(\Gamma _{1}^{\prime })=\mathcal{N}(\Gamma _{1}^{\prime \prime })$. Thus,
the realizing matrix of $\Gamma _{1}^{\prime \prime }$ has required size.
Finally, from Suleimanova criterion is easy to construct the matrices $%
\Gamma _{r_{i}},\Gamma _{1}^{\prime \prime }$ and $\Gamma _{\mathcal{R}},$
and analogously as in Kellogg proof, we use the Brauer Theorem \ref%
{TeoBrauerAA} for construct a nonnegative matrix with spectrum $\Lambda $.
\end{proof}

\begin{remark}
Some realizability criteria, like those of Kellogg and Borobia, give
sufficient conditions for the existence of a nonnegative matrix $A$ with
prescribed spectrum, but not for the construction of $A.$ All realizability
criteria obtained from Brauer or Rado results are constructive, in the sence
that they always allow us to compute a realizing matrix.
\end{remark}

\noindent The following result, due to \v{S}migoc \cite{Smigoc}, can be also
proved as a consequence of Rado's Theorem.
\begin{theorem}
\label{smigoc}\textrm{\cite[Theorem 10]{Smigoc} }\textit{Let }$A=%
\begin{bmatrix}
A_{1} & \mathbf{a} \\
\mathbf{b}^{T} & c%
\end{bmatrix}%
$\textit{\ be an }$n\times n$\textit{\ nonnegative matrix with spectrum }$%
\Lambda $\textit{\ and let }$B$\textit{\ be an }$m\times m$\textit{\
nonnegative matrix with Perron eigenvalue }$\lambda _{1}$\textit{, spectrum }%
$\{\lambda _{1},\Lambda ^{\prime }\}$\textit{\ and maximal diagonal element }%
$d$\textit{. If }$\lambda _{1}\leq c,$\textit{\ then there exists an }$%
(n+m-1)\times (n+m-1)$\textit{\ nonnegative matrix }$M$\textit{\ with
spectrum }$\{\Lambda ,\Lambda ^{\prime }\}$\textit{\ and maximal diagonal
element greater than or equal to }$c+d-\lambda _{1}$\textit{. }
\end{theorem}

\begin{proof}
We suppose without loss of generality that $B\in \mathcal{CS}_{\lambda _{1}}$
with maximal diagonal element $d$ in the position $b_{mm}$. If $\lambda
_{1}\leq c$, we put $\epsilon =c-\lambda _{1}\geq 0$, then from Brauer's
Theorem we see that there exists a nonnegative $B^{\prime }=B+\epsilon
\mathbf{e}\mathbf{e}_{n}$ with spectrum $\{c,\Lambda ^{\prime }\}$ and
maximal diagonal element $d+c-\lambda _{1}$. On the other hand, let $A$ be
nonnegative matrix with spectrum $\Lambda $ and diagonal entries $%
a_{2},a_{3},\ldots ,a_{n},c$. From Rado's Theorem the $(n+m-1)\times (n+m-1)$
matrix
\begin{equation*}
M=%
\begin{bmatrix}
a_{2} &  &  &  &  \\
& a_{3} &  &  &  \\
&  & \ddots  &  &  \\
&  &  & a_{n} &  \\
&  &  &  & B^{\prime }%
\end{bmatrix}%
+XC,\ \ \text{with}\ \ X=%
\begin{bmatrix}
1 & 0 & \cdots  & 0 & 0 \\
0 & 1 & \cdots  & 0 & 0 \\
\vdots  & 0 & \ddots  & \vdots  & \vdots  \\
\vdots  & \vdots  & \ddots  & 1 & 0 \\
\vdots  & \vdots  & \cdots  & 0 & 1 \\
\vdots  & \vdots  & \cdots  & 0 & \vdots  \\
\vdots  & \vdots  & \cdots  & \vdots  & \vdots  \\
0 & \cdots  & \cdots  & 0 & 1%
\end{bmatrix}%
_{(n+m-1)\times n},
\end{equation*}%
has spectrum $\{\Lambda ,\Lambda ^{\prime }\}$ where $\Lambda $ is the
spectrum of $\Omega +CX=A$, being $\Omega =diag\{a_{2},a_{3},\ldots
,a_{n},c\}$. Since $CX=A-\Omega \geq 0,$ it is clear that the $n\times
(n+m-1)$ matrix%
\begin{equation*}
C=\left(
\begin{array}{ccc}
A_{1}^{\prime } & \mathbf{a} & 0\cdots 0 \\
b & 0 & 0\cdots 0%
\end{array}%
\right) ,\ \ (\text{with}\ \ A_{1}^{\prime }\ \ \text{being}\ \ A_{1}\ \
\text{without its diagonal})
\end{equation*}%
is nonnegative. Therefore, the matrix $M$ is nonnegative. Finally, from the
construction of the matrix $M$ it is clear that $M$ has maximal diagonal
element greater than or equal to $c+d-\lambda _{1}$.
\end{proof}

\begin{remark}
\textit{In \cite{Marijuan} the authors construct a map of sufficient
conditions for the RNIEP to have a solution, with inclusion relations or
independency relations between them. There, they point out that Soto 2,
Perfect 2}$^{+}$\textit{, and Soto-Rojo realizability criteria (all them
obtained from Brauer and Rado results) are the most general criteria. In
particular, they conclude that Soto-Rojo criterion contains all
realizability criteria, which are compared in \cite{Marijuan}}\textrm{. }
\end{remark}

\section{SNIEP}

\noindent In this section we consider the symmetric nonnegative inverse
eigenvalue problem (\emph{SNIEP}). It is well known that the \emph{RNIEP}
and the \emph{SNIEP} are equivalent for $n\leq 4,$ while they are different
for $n\geq 5$ \cite{Johnson2}. The first results about symmetric nonnegative
realization are due to Fiedler \cite{Fiedler}. Several realizability
criteria obtained for the \textit{NIEP} have later been proved to be also
symmetric realizability criteria. Fiedler and Radwan, in \cite{Fiedler} and
\cite{Radwan} respectively, show that Kellogg and Borobia realizability
criteria are also symmetric realizability criteria. In \cite{Soto3, Soto5},
the author show that \textit{NIEP} realizability criteria given in \cite%
{Soto} are also symmetric realizability criteria. In \cite{Soto6} the
authors give a symmetric version of the Rado's result, Theorem \ref%
{TeoBrauerD}. Then, by applying Theorem \ref{TeoBrauerD}, they prove a new
symmetric realizability criterion \cite[Theorems $2.6$ and $3.1$]{Soto6},
which strictly contains criteria in \cite{Soto}. In \cite{Soto7} is also
shown that the called Soto p criteria are also \textit{NIEP} and \textit{\
SNIEP} realizability criteria (see \cite{Ellard, Marijuan2}). \newline
\ \newline
Next, by the use of Theorem \ref{TeoBrauerD}, the symmetric Rado version, we
give an alternative proof of the following two results of Fiedler:

\begin{lemma}
\textit{[Fiedler \cite{Fiedler}]\label{FiedLema1} Let }$A$\textit{\ be a
symmetric }$m\times m$\textit{\ matrix with spectrum }$\Lambda _{1}=\{\alpha
_{1},\ldots ,\alpha _{m}\}.$\textit{\ Let }$\mathbf{u=(}u_{1}\ldots ,u_{m})$%
\textit{, }$\left\Vert \mathbf{u}\right\Vert =1,$\textit{\ be a unit
eigenvector of }$A$\textit{\ corresponding to }$\alpha _{1}.$\textit{\ Let }$%
B$\textit{\ be a symmetric }$n\times n$\textit{\ matrix with spectrum }$%
\Lambda _{2}=\{\beta _{1},\ldots ,\beta _{n}\}.$\textit{\ Let }$\mathbf{v=(}%
v_{1},\ldots ,v_{n})$\textit{, }$\left\Vert \mathbf{v}\right\Vert =1,$%
\textit{\ be a unit eigenvector of }$A$\textit{\ corresponding to }$\beta
_{1}.$\textit{\ Then for any scalar }$\rho ,$\textit{\ the matrix }%
\begin{equation*}
C=%
\begin{bmatrix}
A & \rho \mathbf{uv^{T}} \\
\rho \mathbf{vu^{T}} & B%
\end{bmatrix}%
\end{equation*}%
\textit{has spectrum }$\Lambda =\{\gamma _{1},\gamma _{2},\alpha _{2},\ldots
,\alpha _{m},\beta _{2},\ldots ,\beta _{n}\},$\textit{\ where }$\gamma
_{1},\gamma _{2}$\textit{\ are eigenvalues of the matrix }%
\begin{equation*}
\widehat{C}=%
\begin{bmatrix}
\alpha _{1} & \rho  \\
\rho  & \beta _{1}%
\end{bmatrix}.%
\end{equation*}
\end{lemma}

\begin{proof}
The matrix $M=%
\begin{bmatrix}
A & 0 \\
0 & B%
\end{bmatrix}%
$ is symmetric of order $(m+n)$ with eigenvalues $\alpha _{1},\ldots ,\alpha
_{m},\beta _{1},\ldots ,\beta _{n}.$ Let
\begin{eqnarray*}
X_{1}^{T} &=&(u_{1},\ldots ,u_{m},0,\ldots ,0),\text{ and} \\
X_{2}^{T} &=&(0,\ldots ,0,v_{1},\ldots ,v_{n})
\end{eqnarray*}%
$(m+n)-$dimensional vectors. Let $X=[X_{1}\mid X_{2}]$ and $\Omega
=diag\{\alpha _{1},\beta _{1}\}.$ Then $MX=X\Omega $. Let
\begin{equation*}
C=%
\begin{bmatrix}
0 & \rho  \\
\rho  & 0%
\end{bmatrix}.%
\end{equation*}%
Then
\begin{equation*}
XCX^{T}=%
\begin{bmatrix}
0 & \rho \mathbf{uv}^{T} \\
\rho \mathbf{vu}^{T} & 0%
\end{bmatrix}%
\end{equation*}%
and the matrix
\begin{equation*}
M+XCX^{T}=%
\begin{bmatrix}
A & \rho \mathbf{uv}^{T} \\
\rho \mathbf{vu}^{T} & B%
\end{bmatrix}%
\end{equation*}%
is symmetric, and from Theorem \ref{TeoBrauerD}, it has spectrum
\begin{equation*}
\Lambda =\{\gamma _{1},\gamma _{2},\alpha _{2},\ldots ,\alpha _{m},\beta
_{2},\ldots ,\beta _{n}\},
\end{equation*}%
where $\gamma _{1}$ and $\gamma _{2}$ are eigenvalues of the matrix%
\begin{equation*}
\Omega +C=%
\begin{bmatrix}
\alpha _{1} & \rho  \\
\rho  & \beta _{1}%
\end{bmatrix}%
\end{equation*}%
for any $\rho .$
\end{proof}

\ \newline
Observe that Theorem \ref{TeoBrauerD} generalizes Lemma \ref{FiedLema1} of
Fiedler. In fact, if we have symmetric matrices $A_{1},A_{2},\ldots ,A_{p},$
with corresponding spectra $\Lambda _{i}=\{\alpha _{1}^{(i)},\alpha
_{2}^{(i)},\ldots ,\alpha _{ni}^{(i)}\},$ $i=1,2,\ldots ,p,$ and unitary
eigenvectors $\mathbf{u}^{(i)}$ associated, respectively, to the eigenvalues
$\alpha _{1}^{(i)},$ then from Theorem \ref{TeoBrauerD} we may obtain a
symmetric $n\times n$ matrix%
\begin{equation*}
A=(A_{1}\oplus A_{2}\oplus \ldots \oplus A_{p})+XCX^{T},
\end{equation*}%
with spectrum $\{\gamma _{1},\ldots ,\gamma _{p},\alpha _{2}^{(1)},\ldots
,\alpha _{n1}^{(1)},\ldots ,\alpha _{2}^{(p)},\ldots ,\alpha _{np}^{(p)}\},$
where $\gamma _{1},\ldots ,\gamma _{p}$ are eigenvalues of the matrix $%
\Omega +C,$ with $\Omega =diag\{\alpha _{1}^{(1)},\alpha _{1}^{(2)},\ldots
,\alpha _{1}^{(p)}\}.$\newline
\ \newline
In what follows $\mathcal{S}_{n}$ ($\widehat{\mathcal{S}_{n}}$) denote the
set of all lists $\Lambda $ for which there exists an $n\times n$ symmetric
nonnegative (positive) matrix with spectrum $\Lambda $.

\begin{theorem}
\textit{[Fiedler \cite{Fiedler}] If }%
\begin{equation*}
\Lambda _{1}=\{\alpha _{1},\ldots ,\alpha _{m}\}\in \mathcal{S}_{m},\text{ \
\ }\Lambda _{2}=\{\beta _{1},\ldots ,\beta _{n}\}\in \mathcal{S}_{n}
\end{equation*}%
\textit{and }$\alpha _{1}\geq \beta _{1},$\textit{\ then for any }$\epsilon
\geq 0,$\textit{\ }%
\begin{equation*}
\Lambda =\{\alpha _{1}+\epsilon ,\beta _{1}-\epsilon ,\alpha _{2},\ldots
,\alpha _{m},\beta _{2},\ldots ,\beta _{n}\}\in \mathcal{S}_{m+n}.
\end{equation*}
\end{theorem}

\begin{proof}
If $\Lambda _{1}\in \mathcal{S}_{m}$ and $\Lambda _{2}\in \mathcal{S}_{n},$
then there exist symmetric nonnegative matrices $A$ and $B$ with spectrum $%
\Lambda _{1}$ and $\Lambda _{2},$ $A\mathbf{u}=\alpha _{1}\mathbf{u},$ $B%
\mathbf{v}=\beta _{1}\mathbf{v,}$ $\left\Vert \mathbf{u}\right\Vert
=\left\Vert \mathbf{v}\right\Vert =1,$ respectively. Then, as before, from
Theorem \ref{TeoBrauerD} we have
\begin{equation*}
M=%
\begin{bmatrix}
A & 0 \\
0 & B%
\end{bmatrix}%
,\ \ \text{with}\ \ MX=X\Omega ,
\end{equation*}%
where
\begin{equation*}
\Omega =%
\begin{bmatrix}
\alpha _{1} & 0 \\
0 & \beta _{1}%
\end{bmatrix}%
\text{ \ and \ }X=%
\begin{bmatrix}
\mathbf{u} & 0 \\
0 & \mathbf{v}%
\end{bmatrix}%
\end{equation*}%
Moreover, for
\begin{equation*}
C=%
\begin{bmatrix}
0 & \rho  \\
\rho  & 0%
\end{bmatrix}%
,\text{ \ }\rho >0,
\end{equation*}%
\begin{equation*}
M+XCX^{T}=%
\begin{bmatrix}
A & \rho \mathbf{uv^{T}} \\
\rho \mathbf{vu^{T}} & B%
\end{bmatrix}%
\end{equation*}%
is symmetric nonnegative with spectrum $\{\gamma _{1},\gamma _{2},\alpha
_{2},\ldots ,\alpha _{m},\beta _{2},\ldots ,\beta _{n}\}$, where $\gamma _{1}
$ and $\gamma _{2}$ are eigenvalues of $\Omega +C=%
\begin{bmatrix}
\alpha _{1} & \rho  \\
\rho  & \beta _{1}%
\end{bmatrix}%
$. If we choose $\rho =\sqrt{\epsilon (\epsilon +(\alpha _{1}-\beta _{1}))},$
then $\gamma _{1}$ and $\gamma _{2}$ are obtained as $\alpha _{1}+\epsilon $
\ and \ $\beta _{1}-\epsilon ,$ respectively.
\end{proof}

\ \newline
It is known that If $\Lambda =\{\lambda _{1},\lambda _{2},\ldots ,\lambda
_{n}\}\in \mathcal{S}_{n}$ and $\epsilon >0,$ then $\Lambda _{\epsilon
}=\{\lambda _{1}+\epsilon ,\lambda _{2},\ldots ,\lambda _{n}\}\in \widehat{%
\mathcal{S}_{n}}.$ It is clear that this result can be also proved by using
Brauer's result.\newline

\begin{remark}
\textit{The criteria obtained from the symmetric Rado's result are good for
SNIEP (in general Brauer's result destroy the symmetry of a matrix). In \cite%
{Marijuan2} the authors construct a map of sufficient conditions for the
existence of a symmetric nonnegative matrix with prescribed real spectrum.
Again, as for the RNIEP, the most general sufficient conditions for the
SNIEP to have a solution, have been obtained from Brauer or Rado results. In
particular, the criterion given in \cite{Soto6} contains any other
realizability criterion for the SNIEP.}
\end{remark}

\section{Complex NIEP}

In this section we consider the general case in which $\Lambda =\{ \lambda
_{1},\lambda _{2},\ldots ,\lambda _{n}\}$ is a list of complex numbers. In
\cite{Borobia3} the authors give the following complex generalization of
Suleimanova's result. The proof uses Brauer Theorem \ref{TeoBrauerAA}:

\begin{theorem}
\textit{\cite{Borobia3}\label{sucomplex} Let }$\Lambda =\{\lambda
_{1},\lambda _{2},\ldots ,\lambda _{n}\}$\textit{\ with }%
\begin{equation}
\Lambda ^{\prime }=\{\lambda _{2},\ldots ,\lambda _{n}\}\ \subset \ \{\ z\in
\mathbb{C}:Rez\leq 0,\ \left\vert Rez\right\vert \geq \left\vert
Imz\right\vert \}  \label{ComEc1}
\end{equation}%
\textit{Then }$\Lambda $\textit{\ is realizable if and only if }$%
\sum_{i=1}^{n}\lambda _{i}\geq 0.$\textrm{\ }
\end{theorem}

\noindent In \cite{Smigoc}, \v{S}migoc proved that (\ref{ComEc1}) can be
improved to%
\begin{equation}
\Lambda ^{\prime }=\{\lambda _{2},\ldots ,\lambda _{n}\}\ \subset \ \{\ z\in
\mathbb{C}\ :\ Rez\leq 0,\ \sqrt{3}\left\vert Rez\right\vert \geq \left\vert
Imz\right\vert \}.  \label{ComEc2}
\end{equation}%
Then $\Lambda =\{\lambda _{1},\lambda _{2},\ldots ,\lambda _{n}\}$ is
realizable if and only if $\sum_{i=1}^{n}\lambda _{i}\geq 0$.\newline
\ \newline
Next, we give an alternative proof of the \v{S}migoc's result in \cite%
{Smigoc}. First we need the following lemma given in \cite[Theorem 2.2]%
{Soto68}.

\begin{lemma}
\label{Perfect} The numbers $\omega _{1},\omega _{2},\omega _{3}$ and $%
\lambda _{1},\lambda _{2},\lambda _{3}$ $(\lambda _{1}\geq |\lambda
_{i}|,i=2,3)$ are, respectively, the diagonal entries and eigenvalues of a
nonnegative matrix $B\in \mathcal{CS}_{\lambda _{1}}$ if only if \newline
$i)$ $0\leq \omega _{k}\leq \lambda _{1}$, \ \ $k=1,2,3$ \newline
$ii)$ $\omega _{1}+\omega _{2}+\omega _{3}=\lambda _{1}+\lambda _{2}+\lambda
_{3}$ \newline
$iii)$ $\omega _{1}\omega _{2}+\omega _{1}\omega _{3}+\omega _{2}\omega
_{3}\geq \lambda _{1}\lambda _{2}+\lambda _{1}\lambda _{3}+\lambda
_{2}\lambda _{3}$ \newline
$iv)$ $\max \omega _{k}\geq Re\lambda _{2}$.
\end{lemma}

\begin{theorem}
Let $\Lambda =\{\lambda _{1},\lambda _{2},\ldots ,\lambda _{n}\}$, $\Lambda =%
\overline{\Lambda }$, $\lambda _{j}=a_{j}+ib_{j}$ with $a_{j}\leq 0$ , $%
b_{j}>0$ satisfying $b_{j}\leq -\sqrt{3}a_{j}$. Then $\Lambda $ is
realizable if only if $\sum_{k=1}^{n}\lambda _{k}\geq 0$.
\end{theorem}

\begin{proof}
The condition is necessary. Now, suppose that $\sum_{k=1}^{n}\lambda
_{k}\geq 0$. We use induction on n, with $n\geq 2$. For $n=2$, $\Lambda
=\{\lambda _{1},\lambda _{2}\}$ must be a real list with $\lambda _{2}<0$.
Then
\begin{equation*}
A=\frac{1}{2}%
\begin{bmatrix}
\lambda _{1}+\lambda _{2} & \lambda _{1}-\lambda _{2} \\
\lambda _{1}-\lambda _{2} & \lambda _{1}+\lambda _{2}%
\end{bmatrix}%
\end{equation*}%
is nonnegative with spectrum $\Lambda $. For $n=3$ and $\Lambda =\{\lambda
_{1},\lambda _{2},\lambda _{3}\}$ with $\lambda _{j}<0$, $j=2,3$, the
conditions from Lemma \ref{Perfect} are satisfied and therefore there exists
a nonnegative matrix with prescribed eigenvalues and diagonal entries.%
\newline
If $\Lambda =\{\lambda _{1},a+ib,a-ib\}$ with $a<0$, $b\leq -\sqrt{3}a$, $%
\lambda _{1}+2a\geq 0$, then since
\begin{equation*}
2a\lambda _{1}+4a^{2}\leq 0\ \ \text{and}\ \ -3a^{2}+b^{2}\leq 0,
\end{equation*}%
we have $2\lambda _{1}a+a^{2}+b^{2}\leq 0$. So, the conditions from Lemma %
\ref{Perfect} are also satisfied and therefore there exists a $3\times 3$
nonnegative matrix with spectrum $\Lambda $ and the prescribed diagonal
entries. Now, we suppose that lists of Smigoc type, with $m-2$ numbers, $%
4\leq m\leq n$, are realizable. Let
\begin{equation*}
\Lambda ^{\prime }=\{\lambda _{1},\lambda _{2},\ldots ,\lambda _{m}\},\ \
\text{with}\ \ Re\lambda _{j}\leq 0,\ \ \sqrt{3}|Re\lambda _{j}|\geq
|Im\lambda _{j}|,\ \ j=2,\ldots ,m.
\end{equation*}%
We take the partition
\begin{equation*}
\Lambda _{0}=\{\lambda _{1},\lambda _{i},\lambda _{j}\},\ \ \Lambda
_{2}=\Lambda ^{\prime }-\Lambda _{0},\ \ \Lambda _{1}=\Lambda _{3}=\emptyset
,\ \ \text{with}
\end{equation*}%
\begin{equation*}
\Gamma _{2}=\{\lambda _{1}+\lambda _{i}+\lambda _{j}\}\cup \Lambda _{2},\ \
\Gamma _{1}=\Gamma _{3}=\{0\},
\end{equation*}%
where $\lambda _{i},\lambda _{j}$ are real or conjugate complex numbers.
From hypothesis of induction, $\Gamma _{2}$ is realizable by a nonnegative
matrix $A_{2}$. Then
\begin{equation*}
A=%
\begin{bmatrix}
A_{2} &  &  \\
& 0 &  \\
&  & 0%
\end{bmatrix}%
\end{equation*}%
is nonnegative with spectrum $\Gamma _{2}\cup \{0,0\}$. From Lemma \ref%
{Perfect}, we can compute a $3\times 3$ nonnegative matrix $B$ with spectrum
$\Lambda _{0}$ and diagonal entries $\{\lambda _{1}+\lambda _{i}+\lambda
_{j},0,0\}$. Finally, from Rado's Theorem with
\begin{equation*}
X=%
\begin{bmatrix}
\mathbf{x_{2}} & 0 & 0 \\
0 & 1 & 0 \\
0 & 0 & 1%
\end{bmatrix}%
_{m\times 3}\ \text{\ and}\ A_{2}\mathbf{x_{2}}=(\lambda _{1}+\lambda
_{i}+\lambda _{j})\mathbf{x_{2}},
\end{equation*}%
we have $B=\Omega +CX$, $\ \Omega =diag\{\lambda _{1}+\lambda _{i}+\lambda
_{j},0,0\}$, and the $m\times m$ matrix
\begin{equation*}
M=%
\begin{bmatrix}
A_{2} &  &  \\
& 0 &  \\
&  & 0%
\end{bmatrix}%
+XC
\end{equation*}%
has the spectrum $\Lambda =\{\lambda _{1},\ldots ,\lambda _{m}\}$. Moreover
since $A,X$ and $C$ are nonnegative, $M$ is nonnegative.
\end{proof}

\section{Universal Realizability}

\noindent We say that a list of complex numbers $\Lambda =\{\lambda
_{1},\ldots ,\lambda _{n}\}$ is universally realizable (\textit{UR}), if for
each possible Jordan canonical form allowed by $\Lambda $ there is a
nonnegative matrix with spectrum $\Lambda .$ As far as we know, the first
works on the universal realizability problem are due to Minc \cite{Minc2,
Minc3}. In \cite{Minc3} Minc prove that if a list $\Lambda =\{\lambda
_{1},\ldots ,\lambda _{n}\}$ of complex numbers has a positive
diagonalizable realization, then $\Lambda $ is \textit{UR}. Next, we give an
alternative proof for this result of Minc:

\begin{theorem}
Minc \cite{Minc3} \label{Minc2} Let $\Lambda =\{\lambda _{1},\lambda
_{2},\ldots ,\lambda _{n}\}$ be realizable by a positive diagonalizable
matrix $A.$ Then $\Lambda $ is universally realizable.
\end{theorem}

\begin{proof}
Let $A$ be positive with spectrum $\Lambda $ and let $S$ be a nonsingular
matrix such that $S^{-1}AS=J(A)$ is the diagonal Jordan canonical form of $%
A. $ We perturb the diagonal matrix $J(A)$ by using Rado Theorem. It is
clear that the eigenvectors of $J(A)$ are the canonical vectors $\mathbf{e}%
_{1},\mathbf{e}_{2},\ldots ,\mathbf{e}_{n}.$ Let $\Omega =diag\{\lambda
_{2},\lambda _{3},\ldots ,\lambda _{r+1}\}$ and consider the $n\times r$
matrix $X=\left[ \mathbf{e}_{2}\mid \mathbf{e}_{3}\mid \cdots \mid \mathbf{e}%
_{r+1}\right] $ and the $r\times n$ matrix $C$ such that $\Omega +CX$ has
eigenvalues $\lambda _{2},\lambda _{3},\ldots ,\lambda _{r+1}.$ Then $XC$ is
of the form $\sum_{i\in K}E_{i,i+1},$ $K=\{2,3,\ldots ,n-1\},$ and%
\begin{equation*}
J(A)+XC=S^{-1}AS+XC=S^{-1}\left( A+SXCS^{-1}\right) S.
\end{equation*}%
By a convenient ordering of the columns, Minc proved in \cite{Minc3} that
the matrix $SXCS^{-1}$ is real. Hence, for $\epsilon >0$ small enough,
\begin{equation*}
M=A+\epsilon SXCS^{-1}
\end{equation*}%
is positive with Jordan canonical form $J(M)=J(A)+XC.$
\end{proof}

\ \newline
Recently, several results on the universal realizability problem, all which
have been obtained by applying Brauer or Rado results, are given in \cite%
{Soto65, Ccapa1, Ccapa2, Soto8, Diaz, Collao, Julio}. Some of these works
give sufficient conditions for the universal realizability problem for
structured matrices \cite{Soto9, Soto10, Soto11}.


\begin{thebibliography}{99}
\bibitem{Borobia1} A. Borobia, On the Nonnegative Eigenvalue Problem, Linear
Algebra Appl. 223-224 (1995) 131-140.

\bibitem{Borobia3} A. Borobia, J. Moro, R. L. Soto, Negativity compensation
in the nonnegative inverse eigenvalue problem, Linear Algebra Appl 393
(2004) 73-89.

\bibitem{Brauer} A. Brauer, Limits for the characteristic roots of a matrix.
IV: Aplications to stochastic matrices, Duke Math. J., 19 (1952) 75-91.

\bibitem{Ccapa1} J. Ccapa, R. L. Soto, On spectra perturbation and
elementary divisors of positive matrices, Electronic Journal of Linear
Algebra 18 (2009) 462-481.

\bibitem{Ccapa2} J. Ccapa, R. L. Soto, On elementary divisors perturbation
of nonnegative matrices, Linear Algebra Appl 432 (2010) 546-555.

\bibitem{Ciarlet} P. G. Ciarlet, Some results in the theory of nonnegative
matrices, Linear Algebra Appl. 9 (1974) 119-142.

\bibitem{Collao} M. Collao, C. R. Johnson, R. L. Soto, Universal
realizability of spectra with two positive eigenvalues, Linear Algebra Appl.
545 (2018) 226-239.

\bibitem{Diaz} R. C. D\'{\i}az, R. L. Soto, Nonnegative inverse elementary
divisors problem in the left half plane, Linear and Multilinear Algebra 64
(2016) 258-268.

\bibitem{Ellard} R. Ellard, H. \v{S}migoc, Connecting sufficient conditions
for the symmetric nonnegative inverse eigenvalue problem, Linear Algebra
Appl 498 (2016) 521-552.

\bibitem{Fiedler} M. Fiedler, Eigenvalues of nonnegative symmetric matrices,
Linear Algebra Appl. 9 (1974) 119-142.

\bibitem{Guo2} W. Guo, Eigenvalues of nonnegative matrices, Linear Algebra
Appl. 266 (1997) 261-270.

\bibitem{Johnson2} C. R. Johnson, T. J. Laffey, R. Loewy, The real and the
symmetric nonnegative inverse eigenvalue problems are different, Proc. AMS
124 (1996) 3647-3651.

\bibitem{Julio} A. I. Julio, C. Mariju\'{a}n, M. Pisonero, R. L. Soto, On
universal realizability of spectra, Linear Algebra Appl. 563 (2019) 353-372.

\bibitem{Kellogg} R. Kellogg, Matrices similar \ to a positive or
essentially positive matrix, Linear Algebra Appl. 4 (1971) 191-204.

\bibitem{Laffey2} T. J. Laffey, E. Meehan, A characterization of trace zero
nonnegative $5\times 5$ matrices, Linear Algebra Appl. 302-303 (1999)
295-302.

\bibitem{Loewy} R. Loewy, D. London, A note on an inverse problem for
nonnegative matrices, Linear and Multilinear Algebra 6 (1978) 83-90.

\bibitem{Marijuan} C. Mariju\'{a}n, M. Pisonero, R. L. Soto, A map of
sufficient conditions for the real nonnegative inverse eigenvalue problem,
Linear Algebra Appl. 426 (2007) 690-705.

\bibitem{Marijuan2} C. Mariju\'{a}n, M. Pisonero, R. L. Soto, A map of
sufficient conditions for the symmetric nonnegative inverse eigenvalue
problem, Linear Algebra Appl. 530 (2017) 344-365.

\bibitem{Meehan} M. E. Meehan, Some results on matrix spectra. Ph.D. thesis,
National university of Ireland, Dublin, 1998.

\bibitem{Minc2} H. Minc, Inverse elementary divisors problem for doubly
stochastic matrices, Proc. Amer.Math. Soc. 83 (1981) 665-670.

\bibitem{Minc3} H. Minc, Inverse elementary divisors problem for nonnegative
matrices, Linear and Multilinear Algebra 11 (1982) 121-131.

\bibitem{Perfect1} H. Perfect, Methods of constructing certain stochastic
matrices, Duke Math. J. 20 (1953) 395-404.

\bibitem{Perfect2} H. Perfect, Methods of constructing certain stochastic
matrices II, Duke Math. J. 22 (1955) 305-311.

\bibitem{Radwan} N. Radwan, An inverse eigenvalue problem for symmetric and
normal matrices, Linear Algebra Appl. 248 (1996) 101-109.

\bibitem{Reams} R. Reams, An inequality for nonnegative matrices and the
inverse eigenvalue problem, Linear and Multilinear Algebra 41 (1996) 367-375.

\bibitem{Rojo} O. Rojo, R. L. Soto, Existence and construction of
nonnegative matrices with complex spectrum, Linear Algebra Appl. 368 (2003)
53-69.

\bibitem{Salzmann} F. Salzmann, A note on eigenvalues of nonnegative
matrices, Linear Algebra Appl. 5.(1972) 329-338.

\bibitem{Smigoc} H. \v{S}migoc, The inverse eigenvalue problem for
nonnegative matrices, Linear Algebra Appl. 393 (2004) 365-374.

\bibitem{Soto} R. Soto, Existence and construction of nonnegative matrices
with prescribed spectrum, Linear Algebra Appl. 369 (2003) 169-184.

\bibitem{Soto2} R. L. Soto, A. Borobia, J. Moro, On the comparison of some
realizability criteria for the real nonnegative inverse eigenvalue problem,
Linear Algebra Appl. 396.(2005) 223-241.

\bibitem{Soto3} R. L. Soto, Realizability by symmetric nonnegative matrices,
Proyecciones Journal of Mathematics 24 (2005) 65-78.

\bibitem{Soto4} R. L. Soto, O. Rojo, Applications of a Brauer Theorem in the
nonnegative inverse eigenvalue problem, Linear Algebra App. 416 (2006)
844-856.

\bibitem{Soto5} R. L. Soto, Realizability criterion for the symmetric
nonnegative inverse eigenvalue problem, Linear Algebra App. 416 (2006)
783-794.

\bibitem{Soto6} R. L. Soto, O. Rojo, J. Moro, A. Borobia, Symmetric
nonnegative realization of spectra, Electronic Journal of Linear Algebra 16
(2007) 1-18.

\bibitem{Soto65} R. L. Soto, J. Ccapa, Nonnegative matrices with prescribed
elementary divisors, Electronic Journal of Linear Algebra 17 (2008) 287-303.

\bibitem{Soto68} R. L. Soto, M. Salas, C. Manzaneda, Nonnegative realization
of complex spectra, Electronic Journal of Linear Algebra 20 (2010) 595-609.

\bibitem{Soto69} R. L. Soto, O. Rojo, C. B. Manzaneda, On nonnegative
realization of partitioned spectra, Electronic Journal of Linear Algebra 22
(2011) 557-572.

\bibitem{Soto7} R. L. Soto, A family of realizability criteria for the real
and symmetric nonnegative inverse eigenvalue problem, Nunerical Linear
Algebra with Applications 20 (2013) 336-348.

\bibitem{Soto8} R. L. Soto, R. C. D\'{\i}az, H. Nina, M. Salas, Nonnegative
matrices with precribed spectrum and elementary divisors, Linear Algebra
Appl 439 (2013) 3591-3604.

\bibitem{Soto9} R. L. Soto, A. I. Julio, M. Salas, Nonnegative persymmetric
matrices with prescribed elementary divisors, Linear Algebra Appl 483 (2015)
139-157.

\bibitem{Soto10} R. L. Soto, E. Valero, M. Salas, H. Nina, Nonnegative
generalized doubly stochastic matrices with prescribed elementary divisors,
Electronic Journal of Linear Algebra 30 (2015) 704-720.

\bibitem{Soto11} R. L. Soto, R. C. D\'{\i}az, M. Salas, O. Rojo, M-matrices
with prescribed elementary divisors, Inverse Problems 33 (2017) 15 pp.

\bibitem{Spector} O. Spector, A characterization of trace zero symmetric
nonnegative $5\times 5$ matrices, Linear Algebra Appl 434 (2011) 1000-1017.

\bibitem{Suleimanova} H. R. Suleimanova, Stochastic matrices with real
characteristic values, Dokl. Akad. Nauk SSSR 66 (1949) 343-345.

\bibitem{Torre} J. Torre-Mayo, M.R. Abril-Raymundo, E. Alarcia-Est\'{e}vez,
C. Mariju\'{a}n, y M. Pisonero. The nonnegative inverse eigenvalue problem
from the coeffcients of the characteristic polynomial. EBL diagraphs. Linear
Algebra Appl. 426 (2007) 729-773.
\end{thebibliography}
\end{document}